# A $M|M|m|m$ Queue System Transient Behavior Study

Manuel Alberto M. Ferreira

manuel.ferreira@iscte-iul.pt

**Abstract**

It is a very hard task to compute an exact solution for the differential equations, with differences, system that allows the determination of the *M|M|m|m* system transient probabilities. The respective complexity grows with *m*. The computations are extremely fastidious and the length and the fact that the expressions obtained are often approximate, and not exact, will not allow the transient probabilities behavior as time functions characterization. To overcome these problems, in this work it is analyzed how that system can supply approximate values to the *M|M|m|m* queue system. It is also presented an asymptotic method to solve the system that becomes possible in many cases to obtain simple approximated expressions for those probabilities using the *M|M|∞* transient probabilities, very well-known and very much easier to study.

## 1.Introduction

Be a queue system $M|M|m|m$, on what:

- $\lambda_0$ is the Poisson process arrivals rate,

- Each customer receives an exponential distributed service at rate $1/\alpha$,

- There are *m* servers,

- The maximum number of present customers, simultaneously, does not exceed *m*,

- The customers that, when arriving, meet the whole servers occupied are lost for the system.

Call $N(t)$ the number of occupied servers at time $t$ and $P_n(t) = P[N(t) = n], t \geq 0, 0 \leq n \leq m$, the transient probabilities.

To determine the transient probabilities, we are led to the following differential difference equations:

$$\frac{d}{dt}P_n(t) = \lambda_0 P_{n-1}(t) + \frac{1}{\alpha}(n+1)P_{n+1}(t) - \left(\lambda_0 + \frac{1}{\alpha}n\right)P_n(t), 1 \leq n \leq m-1 \quad (1.1),$$

With border conditions:

$$\frac{d}{dt}P_0(t) = \frac{1}{\alpha} P_1(t) - \lambda_0 P_0(t) \quad (1.2),$$

$$\frac{d}{dt}P_m(t) = \lambda_0 P_{m-1}(t) - \frac{1}{\alpha} m P_m(t) \quad (1.3).$$

In this system, the stationary and limit probabilities are identical. They can be obtained from (1.1), (1.2) and (1.3), making $\frac{d}{dt}P_n(t) = 0$, $0 \leq n \leq m$. The solution is:

$$P_n = \frac{\rho_0^n}{n!} \left( \sum_{n=0}^{m} \frac{\rho_0^n}{n!} \right)^{-1}, \quad \rho_0 = \lambda_0 \alpha, \quad n = 0,1,2,\ldots n \quad (1.4).$$

The transient probabilities computation, extremely sensible to the initial values $N(0) = n_0$, demands the resolution of the system integrated by (1.1), (1.2), and (1.3). We can write it in the matrix form, see for instance Dynkin (1982),

$$\bar{P}(t) = e^{tA}\bar{P}(0), t \geq 0 \quad (1.5),$$

Being $A = \begin{bmatrix} -\lambda_0 & \frac{1}{\alpha} & 0 & & & \\ \lambda_0 & -\lambda_0 - \frac{1}{\alpha} & \frac{2}{\alpha} & \cdots & & 0_{3\times 2} \\ 0 & \lambda_0 & -\lambda_0 - \frac{2}{\alpha} & & & \\ & \vdots & & \ddots & & \vdots \\ & & & & -\lambda_0 - \frac{m-1}{\alpha} & \frac{m}{\alpha} \\ & 0_{2\times 3} & & \cdots & \lambda_0 & -\frac{m}{\alpha} \end{bmatrix}$ the

$N(t)$ infinitesimal generator, and $e^{tA} = \sum_{j=0}^{\infty} \frac{t^j}{j!} A^j$ and $\bar{p}(t) = \begin{bmatrix} P_0(t) \\ P_1(t) \\ \vdots \\ P_n(t) \end{bmatrix}$.

For $m = 1$ and $m = 2$ is possible to give (1.5) an alternative form without series:

- **$M|M|1|1$**

$$P_{00}(t) = \frac{\rho_0}{1+\rho_0} e^{-(\lambda_0 + \frac{1}{\alpha})t} + \frac{1}{1+\rho_0}, t \geq 0$$

$$P_{01}(t) = -\frac{\rho_0}{1+\rho_0} e^{-(\lambda_0 + \frac{1}{\alpha})t} + \frac{\rho_0}{1+\rho_0}, t \geq 0$$

$$P_{10}(t) = -\frac{1}{1+\rho_0} e^{-(\lambda_0 + \frac{1}{\alpha})t} + \frac{1}{1+\rho_0}, t \geq 0$$

$$P_{11}(t) = \frac{1}{1+\rho_0} e^{-(\lambda_0 + \frac{1}{\alpha})t} + \frac{\rho_0}{1+\rho_0}, t \geq 0$$

(1.6)

And $P_{ij}(t) = P[N(t) = j | N(0) = i], i, j = 1, 2, \ldots, m, \geq 0$.

- **$M|M|2|2$**

Making $a = \frac{1}{\alpha} + 2\lambda_0, c = \sqrt{\left(\frac{1}{\alpha} - 2\lambda_0\right)^2 + 4\left(2\frac{1}{\alpha^2} + \lambda_0^2\right)}$ and $d = 5\frac{1}{\alpha^2} + 6\lambda_0^2$,

$$P_{00}(t) = \frac{-\lambda_0(ca - d) + \left(\lambda_0^2 + \frac{1}{\alpha}\lambda_0\right)(c - a)}{c(d - a^2)} e^{-\frac{a+c}{2}t}$$

$$+ \frac{-\lambda_0(ca - d) + \left(\lambda_0^2 + \frac{1}{\alpha}\lambda_0\right)(c - a)}{c(d - a^2)} e^{-\frac{a-c}{2}t}$$

$$+ \left(1 + 2\frac{\lambda_0 a - \left(\lambda_0^2 + \frac{1}{\alpha}\lambda_0\right)}{d - a^2}\right), t \geq 0 \quad (1.7)$$

$$P_{01}(t) = \frac{\lambda_0(ca-d) + \left(-2\lambda_0^2 - \frac{1}{\alpha}\lambda_0\right)(c-a)}{c(d-a^2)} e^{-\frac{a+c}{2}t}$$

$$+ \frac{\lambda_0(ca+d) + \left(-2\lambda_0^2 + \frac{1}{\alpha}\lambda_0\right)(c+a)}{c(d-a^2)} e^{-\frac{a-c}{2}t}$$

$$- 2\frac{\lambda_0 a - 2\lambda_0^2 - \frac{1}{\alpha}\lambda_0}{d-a^2}, t \geq 0 \qquad (1.8)$$

$$P_{02}(t) = \lambda_0^2 \frac{c-a}{c(d-a^2)} e^{-\frac{a+c}{2}t} + \lambda_0^2 \frac{c+a}{c(d-a^2)} e^{-\frac{a-c}{2}t} - 2\frac{\lambda_0^2}{d-a^2}, t \geq 0 \quad (1.9)$$

$$P_{10}(t) = \frac{\frac{1}{\alpha}(ca-d) + \left(-2\frac{1}{\alpha}\lambda_0 - \frac{1}{\alpha^2}\right)(c-a)}{c(d-a^2)} e^{-\frac{a+c}{2}t}$$

$$+ \frac{\frac{1}{\alpha}(ca-d) + \left(-2\frac{1}{\alpha}\lambda_0 - \frac{1}{\alpha^2}\right)(c+a)}{c(d-a^2)} e^{-\frac{a-c}{2}t} - \frac{\frac{1}{\alpha}a - 2\frac{1}{\alpha}\lambda_0 - \frac{1}{\alpha^2}}{d-a^2}, t \geq 0 \qquad (1.10)$$

$$P_{11}(t) = \frac{\left(-\lambda_0 - \frac{1}{\alpha}\right)(ca-d) + \left(\frac{5}{\alpha}\lambda_0 + \lambda_0^2 + \frac{1}{\alpha^2}\right)(c-a)}{c(d-a^2)} e^{-\frac{a+c}{2}t}$$

$$+ \frac{\left(-\lambda_0 - \frac{1}{\alpha}\right)(ca+d) + \left(\frac{5}{\alpha}\lambda_0 + \lambda_0^2 + \frac{1}{\alpha^2}\right)(c+a)}{c(d-a^2)} e^{-\frac{a-c}{2}t}$$

$$+ \left(1 + 2\frac{\left(\lambda_0 + \frac{1}{\alpha}\right)a - \left(\frac{5}{\alpha}\lambda_0 + \lambda_0^2 + \frac{1}{\alpha^2}\right)}{d-a^2}\right), t \geq 0 \quad (1.11)$$

$$P_{12}(t) = \frac{\lambda_0(ca-d) + \left(-\lambda_0{}^2 - 3\lambda_0\frac{1}{\alpha}\right)(c-a)}{c(d-a^2)} e^{-\frac{a+c}{2}t}$$

$$+ \frac{\lambda_0(ca+d) + \left(-\lambda_0{}^2 - 3\lambda_0\frac{1}{\alpha}\right)(c+a)}{c(d-a^2)} e^{-\frac{a-c}{2}t}$$

$$- 2\frac{\lambda_0 a - \lambda_0{}^2 - 3\lambda_0\frac{1}{\alpha}}{d-a^2}, t \geq 0 \qquad (1.12)$$

$$P_{20}(t) = 2\frac{1}{\alpha^2}\frac{c-a}{c(d-a^2)} e^{-\frac{a+c}{2}t} + 2\frac{1}{\alpha^2}\frac{c+a}{c(d-a^2)} e^{-\frac{a-c}{2}t} - \frac{\frac{4}{a^2}}{d-a^2},$$

$$t \geq 0 \qquad (1.13)$$

$$P_{21}(t) = \frac{2\frac{1}{\alpha}(ca-d) + \left(-2\frac{1}{\alpha}\lambda_0 - 6\frac{1}{\alpha^2}\right)(c-a)}{c(d-a^2)} e^{-\frac{a+c}{2}t}$$

$$+ \frac{2\frac{1}{\alpha}(ca+d) + \left(-2\frac{1}{\alpha}\lambda_0 - 6\frac{1}{\lambda_0{}^2}\right)(c+a)}{c(d-a^2)} e^{-\frac{a-c}{2}t}$$

$$-2\frac{2\frac{1}{\alpha}a - 2\frac{1}{\alpha}\lambda_0 - \frac{6}{a^2}}{d-a^2}, t \geq 0 \qquad (1.14)$$

$$P_{22}(t) = \frac{-2\frac{1}{\alpha}(ca-d) + \left(2\frac{1}{\alpha}\lambda_0 + 4\frac{1}{\alpha^2}\right)(c-a)}{c(d-a^2)} e^{-\frac{a+c}{2}t}$$

$$+ \frac{-2\frac{1}{\alpha}(ca+d) + \left(2\frac{1}{\alpha}\lambda_0 + 4\frac{1}{\alpha^2}\right)(c+a)}{c(d-a^2)} e^{-\frac{a-c}{2}t} + 1$$

$$-2\frac{-2\frac{1}{\alpha}a + 2\frac{1}{\alpha}\lambda_0 + 4\frac{1}{a^2}}{d-a^2}, t \geq 0 \qquad (1.15)$$

As for $m \geq 3$, it does not seem possible to guarantee obtaining explicit solutions of this type for $P_{ij}(t), i, j = 1, 2, \ldots, m, t \geq 0$. In fact, for this to happen, it is necessary to solve the characteristic equation of A, which is an algebraic equation of degree $m+1$. This equation has only negative solutions and the null solution, corresponding to the stationary solution. Thus, in practice, we only must solve an algebraic equation of degree $m$. Therefore, even from $m = 3$ onwards, we have no guarantee of being able to solve exactly the characteristic equation of A. On the other hand, as shown by the cases in which we obtained exact solutions, the calculations are extremely tedious and the solutions obtained, due to to its extent and approximate value, will not allow the characterization of the behavior of $P_{ij}(t), i, j = 1, 2, \ldots, m, t \geq 0$, as functions of t. They will only be used for approximate calculations.

To solve these problems, let's follow these steps:

-To analyze how (1.5) can provide approximate values for $P_{ij}(t)$, see Ferreira and Ramalhoto (1992).

-To present an asymptotic method for solving the system (1.1), (1.2), and (1.3), due to Knessl (1990), which allows in many cases to obtain simple approximate expressions for the $P_{ij}(t)$ in terms of the transient probabilities of the $M|M|\infty$ system, see Ramalhoto and Ferreira (1991).

-To compare the $M|M|m|m$ and the $M|M|\infty$ systems transient behaviors.

## 1. $M|M|m|m$ Transient Probabilities Matrix Calculation

To establish a numerical method for the $P_{ij}(t)$ calculation, begin to show that if in (1.5) we consider a finite number of terms, F+1, we get:

$$\bar{P}(t) \cong \left(\sum_{j=0}^{F} \frac{t^j}{j!} A^j\right) \bar{P}(0), t \geq 0 \quad (2.1)$$

being the error $R_F(t) = \left(\sum_{j=F+1}^{\infty} \frac{t^j}{j!} A^j\right) \bar{P}(0), t \geq 0$.

Note that:

- (2.1) is an asymptotic development for $\bar{P}(t)$ in the time origin.

neighbor.

- (2.1) can be written as

$$\bar{P}(t) \cong \left( \sum_{j=0}^{F} \frac{(t/\alpha)^j}{j!} B^j \right) \bar{P}(0), t \geq 0 \quad (2.2)$$

where

$$B = \begin{bmatrix} -\rho_0 & 1 & 0 & & & \\ \rho_0 & -\rho_0-1 & 2 & \cdots & & 0_{3\times 2} \\ 0 & \rho_0 & -\rho_0-2 & & & \\ & \vdots & & \ddots & & \vdots \\ & & & & -\rho_0-(m-1) & m \\ & 0_{2\times 3} & & \cdots & \rho_0 & -m \end{bmatrix}.$$

So, with $\rho_0, t$ and $m$ established, (2.1) is asymptotic with $1/\alpha$ in the time origin neighbor.

-In the calculation of (2.1), see for instance Laginha (1958) it is mandatory to carry out

$$\Phi(F, m)) = (F-1)(m+1)^3 + (m+1)^2 + F^2 \quad (2.3)$$

multiplications and

$$\Theta(F, m)) = (F-1)(m+1)^3 + m(m-1) \quad (2.4)$$

sums.

To quantify $R_F$, let's define, see Laginha (1958):

-Norm of a Matrix A, call it $N(A)$, is the greatest sum of the absolute values of the elements of any column of A.

- $N(\lambda A) = |\lambda| N(A), \lambda \in \mathbb{R}$.

- $N(A + B) \leq N(A) + N(B)$.

- $N(AB) \leq N(A)N(B)$.

- $N(A^j) \leq (N(A))^j, j \in \mathbb{N}_0$

So, $N(R_F(t)) = \sum_{j=F+1}^{\infty} \frac{t^j}{j!} (N(A))^j$, since $N(\bar{P}(0)) = 1$ and $t \geq 0$.

Still $N(R_F(t)) \leq \frac{t^F}{F!} (N(A))^F \frac{1}{1-t\frac{N(A)}{F+1}}$, since $F \geq [tN(A)] + 1$ where $[x]$ means $x$ characteristic.

As, in this case, $N(A) \leq 2\left(\lambda_0 + \frac{m}{\alpha}\right)$,

-If $F \geq \left[2t\left(\lambda_0 + \frac{m}{\alpha}\right)\right] + 1$

$$N(R_F(t)) \leq \frac{\left(2t\left(\lambda_0 + \frac{m}{\alpha}\right)\right)^F}{F!} \frac{F+1}{F+1-2t\left(\lambda_0 + \frac{m}{\alpha}\right)} \quad (2.5)$$

The expression (2.5) may be written in the equivalent form:

- If $F \geq \left[\frac{2t}{\alpha}(\rho_0 + m)\right] + 1$

$$N(R_F(t)) \leq \frac{\left(\frac{2t}{\alpha}(\rho_0 + m)\right)^F}{F!} \frac{F+1}{F+1-\frac{2t}{\alpha}(\rho_0 + m)} \quad (2.6)$$

And,

- If $F \geq \left[\frac{2t}{\alpha}m(\rho+1)\right] + 1$

$$N(R_F(t)) \leq \frac{\left(\frac{2t}{\alpha}m(\rho+1)\right)^F}{F!} \frac{F+1}{F+1-\frac{2t}{\alpha}m(\rho+1)}, \rho = \frac{\rho_0}{m} \quad (2.7)$$

Expression (2.5) guarantees that, apart from computational limitations, we can calculate the $P_{ij}(t)$ with the desired approximation if we consider the sufficient terms in (2.1). Ensuring that $N(R_F(t))$ is lesser than a certain $\varepsilon$ value guarantees that each $P_{ij}(t)$ is also calculated with a certain error lesser than $\varepsilon$.

For a given system, that is: fixed $\lambda_0$, $1/\alpha$, and $m$, the situations will be more favorable the smaller the values of $t$ (the fewer terms we will have to consider guaranteeing that the error is lesser than or equal to a certain value). This is in accordance with the fact that (2.1) is asymptotic, at the origin of times, with $t$.

Once $t$, $\rho_0$, and $m$ are fixed, the situations will be more favorable the lower the values of $1/\alpha$. We also saw that, in this case, (2.1) is asymptotic with $1/\alpha$ in the neighborhood of $1/\alpha = 0$.

Finally, situations in which the product $t\,(m/\alpha)(\rho + 1)$ is equal, although they behave identically in relation to the number of terms to be considered to guarantee a certain error upper bound, in fact they are not identical from a computational point of view (because a larger $m$ requires more operations (expressions (2.3) and (2.4)).

In the experiments carried out, to have some sensitivity to these issues, we followed the following method:

-We used (2.5) to determine the first value of $F$ for which $N(R_F(t)) \leq 0.001$,

-For this value of $F$, we calculated the error upper bound given by (2.5) and the number of multiplications and additions given by (2.3) and (2.4), respectively.

In Experiment A,

Table 2.1 **Experiment A**

$M|M|10|10 \quad \lambda_0 = 2 \quad 1/\alpha = 1$

| $t$ | Error U. B. | $F$ | $\Phi$ | $\Theta$ |
|---|---|---|---|---|
| 0.1 | .0005 | 11 | 13552 | 13400 |
| 0.2 | .0004 | 18 | 23072 | 22717 |
| 0.4 | .0005 | 31 | 41012 | 40020 |
| 0.6 | .0005 | 44 | 59290 | 57323 |
| 0.8 | .0005 | 57 | 77906 | 74626 |
| 1.0 | .0005 | 70 | 96860 | 91929 |
| 1.2 | .0005 | 83 | 116152 | 109232 |
| 1.4 | .0005 | 96 | 135782 | 126535 |
| 1.6 | .0005 | 109 | 155750 | 143838 |
| 1.8 | .0005 | 122 | 176056 | 161141 |
| 2.0 | .0005 | 135 | 196700 | 178444 |
| 2.2 | .0005 | 148 | 217682 | 195747 |
| 2.4 | .0005 | 161 | 239002 | 213050 |

it is clear how the number of terms increases with $t$ to the point that from $t = 1.4$ it is necessary to consider $F > 95$. For $t = 2.4$ we have $F = 161$ and it is necessary to carry out around half a million operations.

In Experiments B and C,

## Table 2.2 **Experiment B**

$M|M|10|10 \quad t = 0.1 \quad 1/\alpha = 1$

| $\lambda_0$ | Error U. B. | F | $\Phi$ | $\Theta$ |
|---|---|---|---|---|
| 0.1 | .0004 | 10 | 12200 | 12069 |
| 0.8 | .0008 | 10 | 12200 | 12069 |
| 1.6 | .0003 | 11 | 13552 | 13400 |
| 2.4 | .0007 | 11 | 13552 | 13400 |
| 4.0 | .0006 | 12 | 14906 | 14731 |
| 8.0 | .0009 | 14 | 17620 | 17393 |
| 10.0 | .0003 | 16 | 20342 | 20055 |
| 20.0 | .0006 | 21 | 27182 | 26710 |
| 30.0 | .0003 | 27 | 35456 | 34696 |
| 40.0 | .0005 | 32 | 42406 | 41351 |
| 50.0 | .0009 | 37 | 49406 | 48006 |
| 60.0 | .0005 | 43 | 57872 | 55992 |
| 70.0 | .0008 | 48 | 64982 | 62647 |
| 80.0 | .0004 | 54 | 73580 | 70663 |
| 90.0 | .0006 | 59 | 80800 | 77288 |
| 100.0 | .0010 | 64 | 88070 | 83943 |

## Table 2.3 **Experiment C**

$M|M|20|20 \quad t = 0.1 \quad 1/\alpha = 1$

| $\lambda_0$ | Error U. B. | F | $\Phi$ | $\Theta$ |
|---|---|---|---|---|
| 0.1 | .0003 | 16 | 139612 | 139295 |
| 0.8 | .0005 | 16 | 139612 | 139295 |
| 1.6 | .0009 | 16 | 139612 | 139295 |
| 2.4 | .0004 | 17 | 148906 | 148556 |
| 4.0 | .0004 | 18 | 158202 | 157817 |
| 8.0 | .0005 | 20 | 176800 | 176339 |
| 10.0 | .0006 | 21 | 186102 | 185600 |
| 20.0 | .0003 | 27 | 241956 | 241166 |
| 30.0 | .0005 | 32 | 288556 | 287471 |
| 40.0 | .0009 | 37 | 335206 | 333776 |
| 50.0 | .0005 | 43 | 391252 | 389342 |
| 60.0 | .0005 | 48 | 438012 | 435647 |
| 70.0 | .0004 | 54 | 494190 | 491213 |
| 80.0 | .0006 | 59 | 541060 | 537518 |
| 90.0 | .0010 | 64 | 587980 | 583823 |
| 100.0 | .0005 | 70 | 644350 | 639389 |

it is clear how the increasing of $F$ with $\lambda_0$ in much lesser than the verified in Experiment A with $t$. It is evident that as the greater is $m$ the lesser is the influence of $\lambda_0$ for the same $1/\alpha$.

It should also be noted, comparing Experiments B and C, that although the values of $F$ are not very different, the number of operations is substantially greater in the case of C.

Therefore, it gives the idea that through this process we can obtain good approximations, without great computational effort, for small values, especially of $t$ and $1/\alpha$, but also of $\lambda_0$ and of $m$.

To clarify what is meant by "small values", note that in the Table 2.3 we have:

- $t = 0.1$; $1/\alpha = 1$; $\lambda_0 = 8.0$ and $F = 20$ for an error major of .0005,

- for these values $\frac{tm}{\alpha}(\rho + 1) = 2.8$.

Therefore, for values of $t$, $\alpha$, $m$, and $\rho$ such that $\frac{tm}{\alpha}(\rho + 1) \leq 2.8$ we can guarantee an error of less than .0005, if we consider $F=20$ (a few terms that are not very high compared to those of Experiments A, B, and C). In any case, the hypothesis that (2.1) allows the characterization of $P_{ij}(t)$ as functions of time is ruled out. This further highlights the quality of the $M|M|m|m$ system's approximations by the $M|M|\infty$ system, especially for large values of $m$, which we will see later.

Finally, for the systems $M|M|1|1$ and $M|M|2|2$ we calculated the transient probabilities using the exact expressions, which we determined above, and the computational process that we have been studying (experiments D and E). For the first we took $F=10$ (eleven terms) and for the second $F=15$ (sixteen terms) in (2.1):

Table 2.4 **Experiment D**

$M|M|1|1$ $\quad 1/\alpha = 0.1$

1. $\lambda_0 = 0.005$

| t | $P_{ij}(t)$ | Exact | Approximate |
|---|---|---|---|
| 0.5 | $P_{00}(t)$ | 0.975914 | 0.975914 |
|  | $P_{01}(t)$ | 0.024086 | 0.024086 |
|  | $P_{10}(t)$ | 0.048171 | 0.048171 |
|  | $P_{11}(t)$ | 0.951829 | 0.951829 |
| 1.0 | $P_{00}(t)$ | 0.953569 | 0.953569 |
|  | $P_{01}(t)$ | 0.046431 | 0.046431 |
|  | $P_{10}(t)$ | 0.092861 | 0.092861 |
|  | $P_{11}(t)$ | 0.907139 | 0.907139 |
| 1.5 | $P_{00}(t)$ | 0.932839 | 0.932839 |
|  | $P_{01}(t)$ | 0.067161 | 0.067161 |
|  | $P_{10}(t)$ | 0.134323 | 0.134323 |
|  | $P_{11}(t)$ | 0.865677 | 0.865677 |
| 2.0 | $P_{00}(t)$ | 0.913606 | 0.913606 |
|  | $P_{01}(t)$ | 0.086394 | 0.086394 |

|   |   | $P_{10}(t)$ | 0.172788 | 0.172788 |
|---|---|---|---|---|
|   |   | $P_{11}(t)$ | 0.827212 | 0.827212 |
| 2.5 |   | $P_{00}(t)$ | 0.895763 | 0.895763 |
|   |   | $P_{01}(t)$ | 0.104237 | 0.104237 |
|   |   | $P_{10}(t)$ | 0.208474 | 0.208474 |
|   |   | $P_{11}(t)$ | 0.791526 | 0.791526 |

Equilibrium values: $\begin{cases} P_0 = 0.2 \\ P_1 = 0.8 \end{cases}$

2. $\lambda_0 = 0.4$

| t | $P_{ij}(t)$ | Exact | Approximate |
|---|---|---|---|
| 0.5 | $P_{00}(t)$ | 0.823041 | 0.823041 |
|   | $P_{01}(t)$ | 0.176959 | 0.176959 |
|   | $P_{10}(t)$ | 0.04424 | 0.04424 |
|   | $P_{11}(t)$ | 0.95576 | 0.95576 |
| 1.0 | $P_{00}(t)$ | 0.685225 | 0.685225 |
|   | $P_{01}(t)$ | 0.314775 | 0.314775 |
|   | $P_{10}(t)$ | 0.078694 | 0.078694 |
|   | $P_{11}(t)$ | 0.921306 | 0.921306 |
| 1.5 | $P_{00}(t)$ | 0.577893 | 0.577893 |
|   | $P_{01}(t)$ | 0.422107 | 0.422107 |
|   | $P_{10}(t)$ | 0.105527 | 0.105527 |
|   | $P_{11}(t)$ | 0.894473 | 0.894473 |
| 2.0 | $P_{00}(t)$ | 0.494304 | 0.494304 |
|   | $P_{01}(t)$ | 0.505696 | 0.505696 |
|   | $P_{10}(t)$ | 0.126424 | 0.126424 |
|   | $P_{11}(t)$ | 0.873576 | 0.873576 |
| 2.5 | $P_{00}(t)$ | 0.429204 | 0.429204 |
|   | $P_{01}(t)$ | 0.570796 | 0.570796 |
|   | $P_{10}(t)$ | 0.142699 | 0.142699 |
|   | $P_{11}(t)$ | 0.857301 | 0.857301 |

Equilibrium values: $\begin{cases} P_0 = 0.2 \\ P_1 = 0.8 \end{cases}$

3. $\lambda_0 = 0.9$

| t | $P_{ij}(t)$ | Exact | Approximate |
|---|---|---|---|
| 0.5 | $P_{00}(t)$ | 0.645878 | 0.645878 |
|   | $P_{01}(t)$ | 0.354122 | 0.354122 |
|   | $P_{10}(t)$ | 0.039347 | 0.039347 |
|   | $P_{11}(t)$ | 0.960653 | 0.960653 |
| 1.0 | $P_{00}(t)$ | 0.431091 | 0.431091 |
|   | $P_{01}(t)$ | 0.568909 | 0.568909 |
| | $P_{10}(t)$ | 0.063212 | 0.063212 |
| | $P_{11}(t)$ | 0.936788 | 0.936788 |
| 1.5 | $P_{00}(t)$ | 0.300817 | 0.300819 |
| | $P_{01}(t)$ | 0.699183 | 0.699181 |
| | $P_{10}(t)$ | 0.077687 | 0.077987 |
| | $P_{11}(t)$ | 0.922313 | 0.922313 |
| 2.0 | $P_{00}(t)$ | 0.221802 | 0.221841 |
| | $P_{01}(t)$ | 0.778198 | 0.778159 |
| | $P_{10}(t)$ | 0.086466 | 0.086462 |
| | $P_{11}(t)$ | 0.913534 | 0.913538 |
| 2.5 | $P_{00}(t)$ | 0.173876 | 0.174320 |
| | $P_{01}(t)$ | 0.826124 | 0.825680 |
| | $P_{10}(t)$ | 0.091792 | 0.091742 |
| | $P_{11}(t)$ | 0.908208 | 0.908258 |

Equilibrium values: $\begin{cases} P_0 = 0.1 \\ P_1 = 0.9 \end{cases}$

Table 2.5 **Experiment E**

$M|M|2|2 \quad 1/\alpha = 0.1$

1. $\lambda_0 = 0.4$

| t | $P_{ij}(t)$ | Exact | Approximate |
|---|---|---|---|
| 0.5 | $P_{00}(t)$ | 0.975895 | 0.975910 |
| | $P_{01}(t)$ | 0.023803 | 0.023798 |
| | $P_{02}(t)$ | 0.000302 | 0.000292 |
| | $P_{10}(t)$ | 0.047606 | 0.047596 |
| | $P_{11}(t)$ | 0.929498 | 0.929484 |
| | $P_{12}(t)$ | 0.022896 | 0.02292 |
| | $P_{20}(t)$ | 0.002419 | 0.00234 |
| | $P_{21}(t)$ | 0.091583 | 0.091682 |
| | $P_{22}(t)$ | 0.905998 | 0.905978 |
| 1.0 | $P_{00}(t)$ | 0.953416 | 0.953533 |
| | $P_{01}(t)$ | 0.045412 | 0.045371 |
| | $P_{02}(t)$ | 0.001172 | 0.001096 |
| | $P_{10}(t)$ | 0.090823 | 0.090742 |
| | $P_{11}(t)$ | 0.86728 | 0.867174 |
| | $P_{12}(t)$ | 0.041896 | 0.042084 |
| | $P_{20}(t)$ | 0.009375 | 0.008767 |
| | $P_{21}(t)$ | 0.167584 | 0.168334 |
| | $P_{22}(t)$ | 0.823041 | 0.822899 |
| 1.5 | $P_{00}(t)$ | 0.932344 | 0.932724 |

|   |   | $P_{01}(t)$ | 0.065099 | 0.064964 |
|---|---|---|---|---|
|   |   | $P_{02}(t)$ | 0.002557 | 0.002312 |
|   |   | $P_{10}(t)$ | 0.130198 | 0.129929 |
|   |   | $P_{11}(t)$ | 0.812374 | 0.812042 |
|   |   | $P_{12}(t)$ | 0.057428 | 0.058029 |
|   |   | $P_{20}(t)$ | 0.020456 | 0.018493 |
|   |   | $P_{21}(t)$ | 0.229713 | 0.232118 |
|   |   | $P_{22}(t)$ | 0.749831 | 0.749389 |
| 2.0 |   | $P_{00}(t)$ | 0.912484 | 0.913352 |
|   |   | $P_{01}(t)$ |   |   |
|   |   | $P_{02}(t)$ | 0.083107 | 0.082792 |
|   |   | $P_{10}(t)$ | 0.004413 | 0.003856 |
|   |   | $P_{11}(t)$ | 0.166213 | 0.165584 |
|   |   | $P_{12}(t)$ | 0.763918 | 0.763192 |
|   |   | $P_{20}(t)$ | 0.069868 | 0.071224 |
|   |   | $P_{21}(t)$ | 0.035303 | 0.030849 |
|   |   | $P_{22}(t)$ | 0.279473 | 0.284895 |
|   |   |   | 0.685225 | 0.684256 |
| 2.5 |   | $P_{00}(t)$ | 0.893652 | 0.895298 |
|   |   | $P_{01}(t)$ | 0.099648 | 0.099043 |
|   |   | $P_{02}(t)$ | 0.0067 | 0.005658 |
|   |   | $P_{10}(t)$ | 0.199296 | 0.198087 |
|   |   | $P_{11}(t)$ | 0.721157 | 0.719845 |
|   |   | $P_{12}(t)$ | 0.079547 | 0.082068 |
|   |   | $P_{20}(t)$ | 0.053601 | 0.045266 |
|   |   | $P_{21}(t)$ | 0.31819 | 0.328219 |
|   |   | $P_{22}(t)$ | 0.628209 | 0.62642 |

$$\text{Equilibrium values:} \begin{cases} P_0 = 0.615385 \\ P_1 = 0.307692 \\ P_2 = 0.076923 \end{cases}$$

2. $\lambda_0 = 0.9$

| t | $P_{ij}(t)$ | Exact | Approximate |
|---|---|---|---|
| 0.5 | $P_{00}(t)$ | 0.820932 | 0.822767 |
|   | $P_{01}(t)$ | 0.161724 | 0.160533 |
|   | $P_{02}(t)$ | 0.017344 | 0.016699 |
|   | $P_{10}(t)$ | 0.040431 | 0.040133 |
|   | $P_{11}(t)$ | 0.789173 | 0.790984 |
|   | $P_{12}(t)$ | 0.170396 | 0.168883 |
| | $P_{20}(t)$ | 0.002168 | 0.002087 |
| | $P_{21}(t)$ | 0.085198 | 0.084442 |
| | $P_{22}(t)$ | 0.912634 | 0.913471 |
| 1.0 | $P_{00}(t)$ | 0.670863 | 0.683423 |
| | $P_{01}(t)$ | 0.268106 | 0.260471 |
| | $P_{02}(t)$ | 0.061031 | 0.056106 |
| | $P_{10}(t)$ | 0.067026 | 0.065118 |
| | $P_{11}(t)$ | 0.63452 | 0.646359 |
| | $P_{12}(t)$ | 0.298622 | 0.288524 |
| | $P_{20}(t)$ | 0.007629 | 0.007013 |
| | $P_{21}(t)$ | 0.149311 | 0.144262 |
| | $P_{22}(t)$ | 0.84306 | 0.848725 |
| 1.5 | $P_{00}(t)$ | 0.536257 | 0.572869 |
| | $P_{01}(t)$ | 0.341347 | 0.320478 |
| | $P_{02}(t)$ | 0.122396 | 0.106653 |
| | $P_{10}(t)$ | 0.085337 | 0.080128 |
| | $P_{11}(t)$ | 0.512118 | 0.546076 |
| | $P_{12}(t)$ | 0.402545 | 0.373805 |
| | $P_{20}(t)$ | 0.015299 | 0.013332 |
| | $P_{21}(t)$ | 0,201272 | 0.186902 |
| | $P_{22}(t)$ | 0.783428 | 0.799766 |
| 2.0 | $P_{00}(t)$ | 0.408795 | 0.484422 |
| | $P_{01}(t)$ | 0.394934 | 0.354463 |
| | $P_{02}(t)$ | 0.196271 | 0.161114 |
| | $P_{10}(t)$ | 0.098733 | 0.088616 |
| | $P_{11}(t)$ | 0.408197 | 0.476364 |
| | $P_{12}(t)$ | 0.493069 | 0.435021 |
| | $P_{20}(t)$ | 0.020139 | 0.020139 |
| | $P_{21}(t)$ | 0.246535 | 0.21751 |
| | $P_{22}(t)$ | 0.728932 | 0-76235 |
| 2.5 | $P_{00}(t)$ | 0.283318 | 0.413129 |
| | $P_{01}(t)$ | 0.437069 | 0.371745 |
| | $P_{02}(t)$ | 0.279613 | 0.215126 |
| | $P_{10}(t)$ | 0.109267 | 0.092936 |
| | $P_{11}(t)$ | 0.313857 | 0.427755 |
| | $P_{12}(t)$ | 0.576876 | 0.479308 |
| | $P_{20}(t)$ | 0.034952 | 0.026891 |
| | $P_{21}(t)$ | 0.288438 | 0.239654 |
| | $P_{22}(t)$ | 0.676611 | 0.733455 |

$$\text{Equilibrium values:} \begin{cases} P_0 = 0.076923 \\ P_1 = 0.307692 \\ P_2 = 0.615385 \end{cases}$$

It appears that the values determined by the computational process are quite good. Better for the $M|M|1|1$ system than the $M|M|2|2$ system. In any of the experiments, the values worsen with the increase in $\rho_0$ and $t$.

The calculations presented were carried out using the MATLAB program.

## 3. Asymptotic Approximations for $M|M|m|m$ Transient Probabilities

This section is based on the work of Knessl (1990) that we will present here and on which we will build our results.

To obtain the results of the system integrated by equations (1.1), (1.2), and (1.3), Knessl (1990) began by giving it the form:

$$\epsilon \frac{\partial}{\partial t} P(x,t) = \lambda P(x-\epsilon,t) + \frac{1}{\alpha}(x+\epsilon)P(x+\epsilon,t) - \left(\lambda + \frac{1}{\alpha}x\right)P(x,t)$$

$$\epsilon \frac{\partial}{\partial t} P(0,t) = \epsilon \frac{1}{\alpha} P(\epsilon,t) - \lambda P(0,t) \qquad (3.1)$$

$$\epsilon \frac{\partial}{\partial t} P(1,t) = \lambda P(1-\epsilon,t) - \frac{1}{\alpha} P(0,t)$$

Where $\lambda = \frac{\lambda_0}{m}$, $\epsilon = m^{-1}$ and $x = \epsilon n$ with $P(x,t) = P_n(t)$. Supposed $m$ great, $\lambda_0$ of the same order of magnitude of $m$ ($\lambda_0 = O(m)$) and $\rho = \lambda\alpha < 1$.

Being the problem solution very sensible to initial conditions, supposing that $N(0) = n_0$ with probability 1, Knessl (1990) shows that there are three regions of initial conditions to consider:

1) $n_0 = O(1)$

   The process begins with only a few servers occupied,

2) $n_0 \gg 1$ and $m - n_0 \gg 1$

The process begins with only a strict part of the servers occupied,

3) $m - n_0 = O(1)$

   The process begins with almost the whole servers occupied.

Within each region of initial conditions, the author constructs several asymptotic developments for $P_n(t)$, each being valid in a different part of the plane *(n, t)*. It also calculates an asymptotic development that is valid for small values of *x* and when *N(t)* assumes values close to 0 or $m$. The latter is extremely important to obtain an accurate approximation of the transient blocking probability $P_m(t)$.

Using the "ray method" Knessl (1990) obtains the valid asymptotic development in most of the plane $(n,t)$. The "ray method" was developed by Keller (1978) in the context of geometric optics, but it can also be applied to asymptotic problems of applied probabilities. In this case it consists of assuming that:

$$P(x,t) \sim C(\epsilon) K(x,t) e^{\psi(x,t)/\epsilon} [1 + o(1)], \epsilon \to 0 \quad (3.2).$$

It is verified that $\psi$ satisfies a first-order non-linear partial derivatives equation, solvable by the characteristic method (see Courant and Hilbert (1962)), with the characteristic curves being called "rays". $K(x, t)$ satisfies an equation for linear partial derivatives, and it is necessary, to determine it completely, to force the asymptotic development (3.2) to match that obtained for small values of $x$ and when $N(t)$ is close to $n_0$.

Among the results obtained in this context by Knessl (1990), we highlight Result1 on page 756 and Result 1* on page 757. We emphasize again that it is considered $\lambda_0 = O(m)$, with the developments presented asymptotic when $m \to \infty$.

Let's see then the Result 1. The following notation is used for the system $M|M|m|m$:

$-P_{0n}^m(t) = P[N(t) = n | n_0 = 0], n = 0,1, \ldots, m$,

**Result 1**

*Being $n_0 = 0$,*

   **A.** $n = O(1), t = \epsilon\tau = O(\epsilon)$

$$P_{0n}^m(t) \sim \frac{(\lambda\tau)^n}{n!} e^{-\lambda\tau},$$

   **B.** $n \gg 1, m - n \gg 1$

$$P_{0n}^m(t) \sim \frac{1}{\sqrt{2\pi n}} \exp\left[n \ln\left(\frac{\rho\left(1 - e^{-t/\alpha}\right)}{\epsilon n}\right) + n - \frac{\rho}{\epsilon}\left(1 - e^{-t/\alpha}\right)\right],$$

   **C.** $n = O(1), t = O(1)$

$$P_{0n}^m(t) \sim \frac{\epsilon^{-n}}{n!} \left[\rho\left(1 - e^{-t/\alpha}\right)\right]^n \exp\left[1 - \frac{\rho}{\epsilon}\left(1 - e^{-t/\alpha}\right)\right],$$

**D.** $m - n = l = O(1)$

$$P_{0n}^m(t) \sim \left(\frac{\epsilon}{2\pi}\right)^{1/2} \left\{\left[\rho\left(1 - e^{-t/\alpha}\right)\right]^{-1} - \frac{\rho e^{-t/\alpha}}{\rho - \left(1 - e^{-t/\alpha}\right)^{-1}}\left(1 - e^{-t/\alpha}\right)^l\right\} \cdot \exp\left[\frac{1}{\epsilon}\left(1 - \rho\left(1 - e^{-t/\alpha}\right) + \ln\left[\rho\left(1 - e^{-t/\alpha}\right)\right]\right)\right].$$

If, in Result 1D, $l = 0$ we obtain the asymptotic development for the blocking probability:

$$P_{0m}^m(t) \sim \left(\frac{\epsilon}{2\pi}\right)^{1/2} \frac{\rho\left(1 - e^{-t/\alpha}\right) - \left(1 - e^{-t/\alpha}\right)^{-1}}{\rho - \left(1 - e^{-t/\alpha}\right)^{-1}} \exp\left[\frac{1}{\epsilon}\left(1 - \rho\left(1 - e^{-t/\alpha}\right) + \ln\left[\rho\left(1 - e^{-t/\alpha}\right)\right]\right)\right]. \quad (3.3)$$

For great values of $t$, $1 - e^{-t/\alpha} \cong 1$. And, remembering that $\epsilon = \frac{1}{m}$ and $\rho = \frac{\rho_0}{m}$, we obtain: $P_{0m}^m(t) \cong \frac{1}{\sqrt{2\pi m}} e^m m^{-m} e^{-\rho_0} \rho_0^m$, that is in accordance with the Erlang loss formula in stationary regime:

$$P_m = \frac{\rho_0^m}{m! \sum_{i=0}^{m} \frac{\rho_0^i}{i!}} \quad (3.4)$$

For m great, since, after the Stirling formula $\sqrt{2\pi m} e^{-m} m^m \cong m!$.

Still note that:

-**B.** $n \gg 1, m - n \gg 1$

We also have, for great values of t, that the second member can have the form:

$$\frac{n!}{\sqrt{2\pi n}}\left(\frac{e}{n}\right)^n \frac{e^{-\rho_0}\rho_0^n}{n!} \cong \frac{e^{-\rho_0}\rho_0^n}{n!},$$ for n great enough, that is the $M|G|\infty$ system equilibrium distribution, see Ferreira and Andrade (2009).

As for the Result 1*, where Knessl !990) considers $n_0 \geq 0$ and $n_0 = O(1)$, we are going to present it concretized for $n_0 = 1$. Designating it, in this version, by Result 2 and putting for the $M|M|m|m$ system:

$-P_{1n}^m(t) = P[N(t) = n|n_0 = 1], n = 0,1,\ldots,m,$

**Result 2**

*Being $n_0 = 1$,*

A. $n = O(1), \ t = O(m^{-1})$

$$P_{1n}^m(t) \sim e^{-\lambda_0 t}\frac{(\lambda_0 t)^{n-1}}{(n-1)!}, n \geq 1,$$

$$P_{1n}^m(t) \sim e^{-\lambda t}\left(\frac{t}{\alpha}\right)^{1-n}, n \leq 1,$$

B. $n \gg 1, m - n \gg 1$

$$P_{1n}^m(t) \sim \frac{1}{\sqrt{2\pi n}}\left[1 - e^{-t/\alpha} + \frac{n}{\rho_0(e^{t/\alpha}-1)}\right]\left(\frac{e}{n}\right)^n e^{-\rho_0\left(1-e^{-t/\alpha}\right)}\left(\rho_0\left(1-e^{-t/\alpha}\right)\right)^n,$$

C. $n = O(1), t = O(m^{-1/2})$

$$P_{1n}^m(t) \sim e^{-m\lambda_0 t + \frac{1}{2}\frac{\lambda_0}{\alpha}t^2}\frac{1}{\sqrt{m^{1-n}}}\sum_{l=0}^{n} n \left(\frac{\rho_0}{m}\right)^{n-l}\frac{\left(\frac{t}{\alpha}\sqrt{m}\right)^{n+1-2l}}{(n-1)!},$$

D. $n = O(1), t = O(1)$

$$P_{1n}^m(t) \sim \left(1 - e^{-t/\alpha}\right) \frac{\left[\rho_0\left(1 - e^{-t/\alpha}\right)\right]^n}{n!} e^{-\rho_0\left(1 - e^{-t/\alpha}\right)},$$

**E.** $m - n = l = O(1)$

$$P_{1n}^m(t) \sim \frac{1}{\sqrt{2\pi m}} \left(1 - e^{-t/\alpha} + \frac{m}{\rho_0(e^{t/\alpha} - 1)}\right) \cdot \left[\left(\frac{\rho_0\left(1 - e^{-t/\alpha}\right)}{m}\right)^{-l}\right.$$

$$\left. - \frac{\rho_0 e^{-t/\alpha}}{\rho_0 - m\left(1 - e^{-t/\alpha}\right)^{-1}} \left(1 - e^{-t/\alpha}\right)^l \right] \left(\frac{e}{m}\right)^m \left[\rho_0\left(1 - e^{-t/\alpha}\right)\right]^m e^{-\rho_0\left(1 - e^{-t/\alpha}\right)}$$

If, in Result 2 E, $l = 0$ we obtain the asymptotic development for the blocking probability:

$$P_{1m}^m(t) \sim \frac{1}{\sqrt{2\pi m}} \frac{\rho_0\left(1 - e^{-t/\alpha}\right) - m\left(1 - e^{-t/\alpha}\right)^{-1}}{\rho_0 - m\left(1 - e^{-t/\alpha}\right)^{-1}} \left(\frac{e}{m}\right)^m \left(1 - e^{-t/\alpha}\right.$$

$$\left. + \frac{m}{\rho_0(e^{t/\alpha} - 1)}\right) e^{-\rho_0\left(1 - e^{-t/\alpha}\right)} \left[\rho_0\left(1 - e^{-t/\alpha}\right)\right]^m \quad (3,5).$$

For large values of $t$, $1 - e^{-t/\alpha} \cong 1$ and $\frac{m}{\rho_0(e^{t/\alpha}-1)} \cong 0$. Thus, the second member of (3.5), for large t is approximately equal to $\frac{1}{\sqrt{2\pi m}} \left(\frac{e}{m}\right)^m e^{-\rho_0} \rho_0^m \cong \frac{e^{-\rho_0} \rho_0^m}{m!}$, which is also in accordance with the Erlang loss formula, for large m in steady state.

Still note that:

**-B.** $n \gg 1, m - n \gg 1$

We also have, for great values of t, that the second member can have the form:
$\frac{n!}{\sqrt{2\pi n}} \left(\frac{e}{n}\right)^n \frac{e^{-\rho_0} \rho_0^n}{n!} \cong \frac{e^{-\rho_0} \rho_0^n}{n!}$, for n great enough, that is the $M|G|\infty$ system equilibrium distribution, see Ferreira and Andrade (2009).

## 4. $M|M|m|m$ Transient Probabilities contrast with $M|M|\infty$ Transient Probabilities

Consider a $M|M|\infty$ queue system where:

-$\lambda_0$ is the arrivals Poisson process rate,

-Each customer receives an exponential service time with rate $1/\alpha$.

After Ferreira () we have for the $M|M|\infty$ queue:

$$P_{0n}^{\infty}(t) = \frac{e^{-\rho_0\left(1-e^{-t/\alpha}\right)}}{n!}\left[\rho_0\left(1-e^{-t/\alpha}\right)\right]^n, t \geq 0, n = 0,1,2,\ldots (4.1),$$

$$P_{1n}^{\infty}(t) = \left(1 - e^{-t/\alpha} + \frac{n}{\rho_0\left(e^{t/\alpha}-1\right)}\right)\frac{e^{-\rho_0\left(1-e^{-t/\alpha}\right)}}{n!}\left[\rho_0\left(1-e^{-t/\alpha}\right)\right]^n, t \geq 0, n = 0,1,2,\ldots (4.2),$$

Em que $\rho_0 = \lambda_0\alpha$.

So, from (4.1) and Result 1, we can conclude immediately that:

**Theorem 4.1**

For a $M|M|m|m$ queue system on what $\rho < 1$ and $\lambda_0 = O(m)$

A. $n \gg 1, m - n \gg 1$

$$P_{0n}^m(t) \sim \frac{n!}{\sqrt{2\pi n}}\left(\frac{e}{n}\right)^n P_{0n}^{\infty}(t), m \to \infty,$$

B. $n = O(1), t = O(1)$

$$P_{0m}^m(t) \sim P_{0n}^{\infty}(t), m \to \infty,$$

C. Blocking Probability

$$P_{0m}^m(t) \sim K(t,m) P_{0m}^{\infty}(t), m \to \infty, K(t,m) \leq \frac{m!}{\sqrt{2\pi m}}\left(\frac{e}{m}\right)^m.$$

And, from (4.2) and Result 2, we can conclude immediately that:

**Theorem 4.2**

For a $M|M|m|m$ queue system on what $\rho < 1$ and $\lambda_0 = O(m)$

    **A.** $n \gg 1, m - n \gg 1$

$$P_{1n}^m(t) \sim \frac{n!}{\sqrt{2\pi n}} \left(\frac{e}{n}\right)^n P_{1n}^\infty(t), m \to \infty,$$

    **B.** $n = O(1), t = O(1)$

$$P_{1n}^m(t) \sim K(t,n) P_{1n}^\infty(t), m \to \infty, K(t,n) \leq 1),$$

    **C.** Blocking Probability

$$P_{1m}^m(t) \sim K(t,m) P_{1m}^\infty(t), m \to \infty, K(t,m) \leq \frac{m!}{\sqrt{2\pi m}} \left(\frac{e}{m}\right)^m.$$

These theorems show, beyond any doubt, that for large $m$ there is an evident relationship between the transient behavior of the $M|M|m|m$ and $M|M|\infty$ systems. Like this,

-For $n \gg 1$, $m-n \gg 1$, that is: the system has a strict part of busy servers (cases **A.**), $P_{0n}^m(t)$ and $P_{1n}^m$ are well approximated by expressions that, apart from a multiplicative factor, which depends only on $n$, are identical to $P_{0n}^\infty(t)$ and $P_{1n}^\infty$, respectively, for $m$ great,

-For $(n, t)$ confined to a neighborhood of $(0,0)$ (**B.** cases), $P_{0n}^\infty(t)$ is a good approximation of $P_{0n}^m(t)$ such as $P_{1n}^\infty(t)$ is a good approximation of $P_{1n}^m(t)$, for $m$ great,

-The blocking probabilities $P_{0m}^m(t)$ and $P_{1m}^m(t)$ are well approximated by expressions that are always inferior or equal to others that, apart from a multiplicative factor that depends only on $m$, are identical to $P_{0m}^\infty(t)$ and $P_{1m}^\infty(t)$ respectively (cases **C.**), for $m$ great.

Therefore, in these situations, the probabilities $P_{0n}^m(t)$ and $P_{1n}^m(t)$ are well approximated by expressions that are identical or inferior or equal to others that have behaviors, such as functions of $t$, identical to those of $P_{0n}^\infty(t)$ and $P_{1n}^\infty(t)$, respectively.

In this way, behaviors as functions of $t$ of $P_{0n}^\infty(t)$ and $P_{1n}^\infty$ give indications about those of $P_{0n}^m(t)$ and of $P_{1n}^m(t)$ respectively, under the indicated conditions. And the expressions for $P_{0n}^m(t)$ and $P_{1n}^m(t)$ are much simpler than the expressions obtained for the solutions of the system of equations given by (1.1), (1.2), and (1.3).

## 5.Conclusions

Due to the fact that the solutions of the system constituted by equations (1.1), (1.2), and (1.3), for $m \geq 3$, are difficult to obtain analytically, if not impossible, in addition to the fact that, when possible, the expressions obtained are difficult to analyze in terms of their behavior as functions of t, we tried to obtain these solutions using numerical methods.

The matrix method, initially discussed, does not present major problems from the point of view of mere calculation, since we can calculate the values and determine upper bounds of the respective errors. Its only possible problem is the computational effort it requires, but given the current capabilities available, we believe it can be easily overcome. But the behavior of expressions as functions of t can only be intuited through tedious analysis of series of numbers.

Thus, the results of Theorem 9.1 and Theorem 9.2 are attractive because they present simple formulas to approximately calculate the $P_{0n}^m(t)$ and $P_{1n}^m(t)$, respectively and allow the relationship of the transient behaviors of the systems $M|M|m|m$ and $M|M|\infty$. They have the flaw that it is not possible to determine upper bounds of the errors made when we use them in the approximate calculation of $P_{0n}^m(t)$ or $P_{1n}^m(t)$.

## 6.References

bibliography**1**.Dynkin, E. B., "Markov processes and related problems of analysis (related papers)". London: Cambridge University Press. 1982.

**2**.Ferreira, M. A. M. and Ramalhoto, M. F., "Cálculo Matricial das Probabilidades Transeuntes do Sistema $M|M|m|m$". Comunicação apresentada no V Congresso da APDIO, Évora. 1992.

**3**.Knessl, C., "On the Transient Behavior of the $M|M|m|m$ Loss Model". Communications in Statistics – Stochastic Models, 6(4), 749-776, 1990.

**4**.Ramalhoto, M. F. and Ferreira, M. A. M., "Comportamento Transeunte de Sistemas Markovianos de Perda Pura". Actas da 3ª Conferência sobre Aplicação da Matemática à Economia e à Gestão. ISEG, Universidade Técnica de Lisboa. 1991.

**5**.Laginha, J. J., "A resolução do sistema de Leontieff". Separata dos Anais do ISCEF-Tomo1-Volume XXVI, 1958.

**6**. Keller, J. B., "Rays, waves and asymptotics". Bull. Am. Math. Soc., 84, 727-750, 1978.

**7**. Courant, R. and Hilbert, D., "Methods of Mathematical Physics, vol. 2". Interscience. New York. 1962.

**8**. Ferreira, M. A. M. and Andrade, M. A. P., "The ties between the M/G/oo queue system transient behaviour and the busy period". IJAR-International Journal of Academic Research, 1(1), 84-92, 2009.